%
%
%
\input amstex.tex
\input amsppt.sty
\input epsf

\NoBlackBoxes
\advance\voffset+1cm
\advance\hoffset+1cm

\topmatter
\title A short proof of the Twelve points theorem \endtitle
			     \author
			     M. Cencelj, D. Repov\v s and M. Skopenkov
			     \endauthor
			     \address
			     Institute for Mathematics, Physics and Mechanics, University of
			     Ljubljana, P. O. Box 2964, 1001 Ljubljana,
			     Slovenia.  e-mail: matija.cencelj\@uni-lj.si \endaddress
			     \address Institute for Mathematics, Physics and Mechanics, University of
			     Ljubljana, P. O. Box 2964, 1001 Ljubljana,
			     Slovenia.  e-mail: dusan.repovs\@uni-lj.si \endaddress
			     \address
			     Department of Differential Geometry, Faculty of Mechanics and Mathematics,
			     Moscow State University, Moscow, Russia 119992.
e-mail: skopenkov\@rambler.ru \endaddress
\author M. Cencelj, D. Repov\v s and M. Skopenkov  \endauthor
\abstract We present a short elementary proof of
the following  Twelve Points Theorem:
Let $M$ be a convex polygon with vertices at the lattice points,
containing a single lattice point in its interior.
Denote by $m$ (resp. $m^*$)
the number of lattice points in the boundary of $M$
(resp. in the boundary of the dual polygon).
Then
$$
m+m^*=12.
$$
\endabstract
\subjclass 52B20 \endsubjclass
\keywords Lattice, lattice polygon, dual polygon, the Pick formula,
toric varietes
\endkeywords
\thanks
Repov\v s and Cencelj were supported in part by the Ministry for Education,
Science and Sport of the Republic of Slovenia Research Program No.~101-509.
Skopenkov was supported in part by Moscow grant and
grant ''Geometric topology'' No.~02-01-00014
of Russian Foundation of Basic Research.
\endthanks
\endtopmatter

\document

The Twelve points theorem is an elegant theorem, which is easy to formulate, but
no simple proof was available until now.
In this paper we present a short and elementary proof of this result.
To state our theorem we need the following

\definition{Definition of the dual polygon}
Let $M=A_1A_2\dots A_n$ be a convex polygon all of whose vertices
lie in the lattice of points with integer coordinates (Figure~1 on the left).
Suppose that $O$ is the only lattice point in the interior of $M$. Draw
the vectors
$\overarrow {A_1A_2}, \overarrow {A_2A_3}, \dots \overarrow {A_nA_1}$
from the point $O$. Choose on each of the obtained segments
the nearest to $O$ lattice point distinct from $O$.
Connecting the $n$ chosen points consequtively, we get a
polygon $M^*$, {\it dual} to the original polygon
(Figure~1 on the right).
Denote by $m$ the number of lattice points in the boundary of $M$,
and by $m^*$ --- in the boundary of the dual polygon.
\enddefinition

\bigskip

$$ \matrix
\epsfysize=3cm
\epsfbox{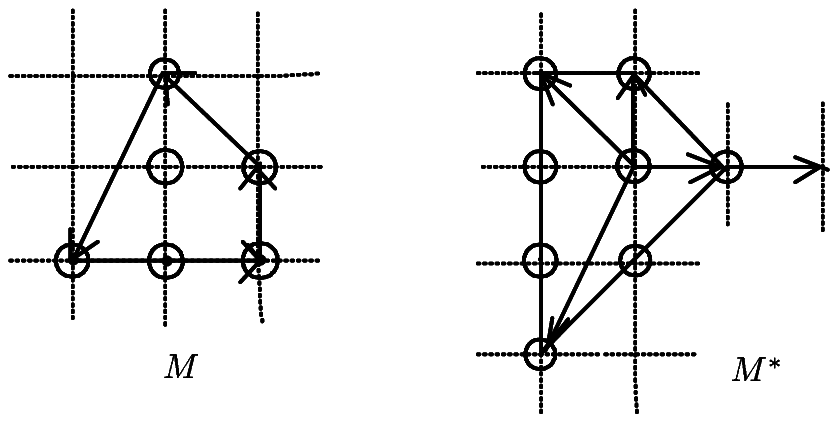}\\
\endmatrix
$$

\nopagebreak
\smallskip

\centerline{Figure~1.}

\bigskip

\proclaim{The Twelve Points Theorem}
Suppose that $O$ is the only lattice point in the interior of a convex
polygon $M$; then
$$
m+m^*=12.
$$
\endproclaim

This theorem appeared in \cite{Ful93}.
There are only some hints for the proof, applying the theory of toric varietes.
In an interesting paper \cite{PRV00}, completely dedicated to the 12 points
theorem,
even four different  proofs are discussed. Two of them are rather long and
they use toric varieties and modular forms respectively. There are
also outlined two proofs applying only linear algebra.
The first of them is exhausting (there are 16 different types
of polygons $M$ in our theorem up to
$\operatornamewithlimits{SL}_2(\Bbb Z)$).
The idea of the second one is very close to the proof of recent paper.

Our elementary proof is analogous to one of the proofs of the Pick formula.
We reduce the Twelve points theorem to the specific case when
$M$ is a parallelogram and $m=4$. Let us begin with this latter case.
\smallskip
(1) {\it If $M=ABCD$ ia a parallelogram without lattice points in its sides
then $m+m^*=12$ (Figure~2).}
\smallskip
Indeed, in this case $O=AC\cap BD$ because the point
symmetric to the point $O$ with respect to $AC\cap BD$ is a lattice point and
belongs to the interior of $ABCD$, so it coincides with $O$. It is easy to
show that $M^*$ is a parallelogram with sides obtained from the diagonals
$AC$ and $BD$ by parallel translations with vectors
$\pm \overrightarrow{OB}$  nd $\pm \overrightarrow{OA}$, respectively.
Since a unique lattice point $O$ belongs to these diagonals, then any side of
the parallelogram $M^*$ contains one lattice point, hence, $b+b^*=4+8=12$.

\bigskip

$$\matrix
\epsfysize=3cm
\epsfbox{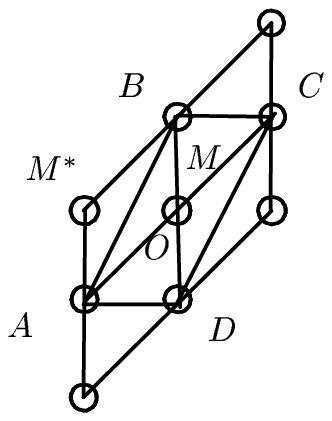}
\endmatrix
$$

\nopagebreak
\smallskip

\centerline{Figure~2.}

\bigskip

Now suppose that $M=A_1A_2\dots A_n$.
Let us assume that all the lattice
boundary points of $M$ are vertices (possibly, with the angle $180^\circ$).
This does not affect the definition of $M^*$.
Assume that some triangle $A_{i-1}A_iA_{i+1}$
is {\it simple}, i. e. it contains no lattice points except its vertices
(neither in the interior nor in the boundary).
{\it Deleting a triangle} is cutting off from polygon $M$ the triangle
$A_{i-1}A_iA_{i+1}$. The reverse operation is called {\it adding a triangle}.
Our reduction is based on the following assertion:
\smallskip
(2) {\it The value $m+m^*$ is preserved under delleting or adding a triangle.}
\smallskip
It is sufficient to prove that deleting a simple triangle,
say $A_1A_2A_3$, from $M$ gives adding a simple triangle
$A_{12}A_{13}A_{23}$ to $M^*$ (Figure~3).
Here by $A_{kl}$ we denote the point such that
$\overrightarrow{OA_{kl}}=\overrightarrow{A_kA_l}$.
In particular, if $l=k+1$ then $A_{kl}$ is a vertex of the polygon $M^*$.
Delete $A_1A_2A_3$. Then one should delete from $M^*$ the vertices
$A_{12}$ and $A_{23}$, and add to it a new vertex $A_{13}$.
The last vertex shold be joined by segments with $A_{n1}$ and $A_{34}$.

Let us show that the points $A_{12}$ and $A_{23}$ beling to these segments.
Indeed, since $O$ is the only lattice point inside $M$,
it follows that the triangles $A_1OA_3$, $A_2OA_3$, $A_4OA_3$ are simple.
By the Pick formula their areas are equal to $1/2$.
Since they have a common base $OA_3$, it follows that the projections
of the vectors $\overrightarrow{A_1A_3}$,
$\overrightarrow{A_1A_3}$ and $\overrightarrow{A_1A_3}$
on the direction normal to $OA_3$ are equal.

This implies that the points $A_{13}$, $A_{23}$ and $A_{34}$
belong to the same line, and $A_{23}$ lies between the two others,
because $M$ is convex.
It can be proved analogously that $A_{12}$ belongs to the segment
$A_{n1}A_{13}$.
Therefore the transformation of $M^*$ is just adding
the triangle $A_{12}A_{13}A_{23}$.

Now note that the triangle $OA_{12}A_{13}$
is obtained from a simple triangle $A_1A_2A_3$
by a parallel translation, and $OA_{23}A_{13}$ is obtained from it
by a cental symmetry. So the triangle
$A_{12}A_{13}A_{23}$ is simple, and assertion (2) is proved.

\bigskip

$$\matrix
\epsfysize=3cm
\epsfbox{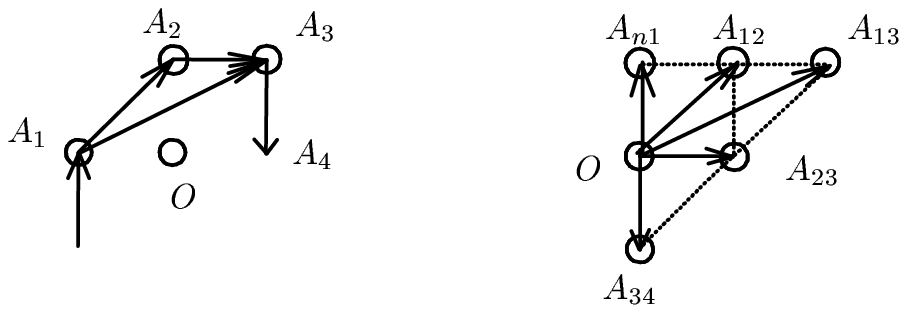}
\endmatrix
$$

\nopagebreak
\smallskip

\centerline{Figure~3.}

\bigskip

For the proof of our theorem it remains to notice the following:
\smallskip
(3) {\it From any polygon $M$ one
can obtain a parallelogram without lattice points in the sides
by a sequence of deleting and adding triangles.}
\smallskip
Indeed, first assume that $M$ has a diagonal
not passing through $O$.
Cut $M$ along this diagonal and consider the obtained part
not containing $O$.

This part necessarily contains a simple triangle of the form
$A_{i-1}A_iA_{i+1}$.
Deleting it we decrease the number $m$.
Repeat this operation until it is possible.
Repetition is impossible only in the following 3 cases
(when such a diagonal does not exist):

\flushpar A) $m=4$, $M=ABCD$, $O=AC\cap BD$.
Since the segments $OA$, $OB$, $OC$ and $OD$
do not contain lattice points, then $OA=OC$ and $OB=OD$,
that is $ABCD$ is the required parallelogram.

\flushpar B) $m=4$, $M=ABD$, and $C$ belongs to the segment $BD$.
In this case let us denote by $D'$
the point symmetric to $D$ with respect to $O$,
and denote by $E$ the middle point of $D'B$.
The required sequence of deleting/adding of triangles has the form:
$$
ABCD\to AEBCD\to AD'EBCD\to AD'ECD\to AD'CD \quad\text{(Figure~4 on the left).}
$$

\flushpar C) $m=3$, $M=ABC$.
In this case denote by $A'$ and $C'$
the points symmetric to $A$ and $C$, respectively, with respect to $O$.
The required sequence of deleting/adding of triangles has the form:
$$
ABC\to AC'BC\to AC'BA'C\to AC'A'C \quad\text{(Figure~4 on the right).}
$$

So in each case we obtain the required parallelogram,
that completes the proof of our theorem.

\bigskip

$$\matrix
\epsfysize=3cm
\epsfbox{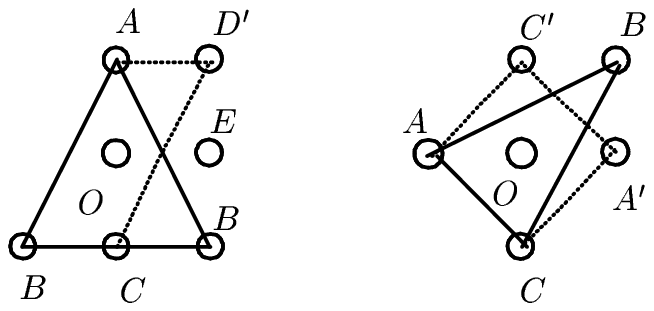}
\endmatrix
$$

\nopagebreak
\smallskip

\centerline{Figure~4.}

\bigskip

{\bf Aknowledgement.} The authors are grateful to V.V. Prasolov
for the statement of the problem.

\enddocument
\end


\author{ УДК 514 \\
Д. Реповш, М. Скопенков и М. Ценцель}
\title{Элементарное доказательство теоремы
о 12 целых точках}
\address{Institute for Mathematics, Physics and Mechanics, University of
Ljubljana, P. O. Box 2964, 1001 Ljubljana,
Slovenia. E-mail: dusan.repovs@uni-lj.si.}
\address{Москва 119992, Московский Государственный Университет,
механико-математический факультет, кафедра дифференциальной
геометрии и приложений. E-mail: skopenkov@rambler.ru.}
\address{Institute for Mathematics, Physics and Mechanics, University of
Ljubljana, P. O. Box 2964, 1001 Ljubljana,
Slovenia. E-mail: matija.cencelj@uni-lj.si.}
\thanks{Реповш и Ценцель частично поддержаны
Министерством Образования, Науки и Спорта Республики Словения,
исследовательская программа No.~101-509.
Скопенков частично поддержан грантом РФФИ No.~02-01-00014.}
\subjclass[2000]{52B20}
\maketitle

Теорема о 12 целых точках --- красивое и простое по формулировке
утверждение, у которого до последнего времени не было простого
доказательства. Данная заметка посвящена ее элементарному доказательству.
Для формулировки нам понадобится

\begin{definition*}[Определение двойственного многоугольника]
Пусть на плоскости фиксирована декартова система координат.
Рассмотрим выпуклый многоугольник $M=A_1A_2\dots A_n$
с вершинами в целых точках, содержащий внутри себя
ровно одну целую точку $O$ (Рис.~1).
Отложим от точки $O$ вектора сторон многоугольника
$\overrightarrow{A_1A_2},\overrightarrow{A_2A_3},
\dots,\overrightarrow{A_nA_1}$ и на каждом из
них отметим ближайшую к $O$ целую точку (отличную от $O$).
Соединяя последовательно отмеченные точки,
получим многоугольник $M^*$, {\it двойственный} исходному.
Обозначим количество целых точек на границах многоугольников $M$ и $M^*$
через $m$ и $m^*$ соответственно.
\end{definition*}

\bigskip

$$
\epsfysize=3cm
\epsfbox{1.eps}
$$

\nopagebreak
\smallskip

\centerline{Рисунок~1.}

\bigskip


\begin{theorem}[о 12 целых точках]
Пусть $M$ --- выпуклый многоугольник с вершинами в целых точках,
содержащий внутри ровно одну целую точку.
Тогда
$$
m+m^*=12.
$$
\end{theorem}

Эта теорема была сформулирована в книге \cite{Ful93}
с некоторыми указаниями к доказательству,
использующему технику торических многообразий (см. также \cite{Kho97}).
В интересной статье \cite{PRV00}, целиком посвященной
этой теореме, обсуждаются 4 различных доказательства.
Три из них довольно объемны и основаны соответственно
на сведении к перебору 16 вариантов, теории модулярных форм
и теории торических многообразий.
Идея четвертого очень близка к идее данной заметки.
(Для него в указанной статье приводится только план, прямая реализация
которого также довольно громоздка.)

Приведенное ниже элементарное доказательство аналогично
одному из доказательств формулы Пика. Мы сводим теорему
о 12 целых точках к частному случаю, когда $M$ --- параллелограмм и
$m=4$. Начнем именно с этого случая:

\smallskip
(1) {\it Если $M=ABCD$ --- параллелограмм без целых точек на сторонах,
то $m+m^*=12$ (Рис.~2 ).}
\smallskip

Действительно, в этом случае $O=AC\cap BD$. Это следует из того, что
точка, симметричная точке $O$ относительно $AC\cap BD$, ---
целая и лежит внутри $ABCD$, а поэтому совпадает с $O$.
Легко видеть, что $M^*$ --- параллелограмм, стороны которого
получаются из диагоналей $AC$ и $BD$ параллельными переносами
на вектора $\pm \overrightarrow{OB}$ и $\pm \overrightarrow{OA}$
соответственно.
Так как на этих диагоналях лежит единственная целая точка $O$,
то на каждой стороне параллелограмма $M^*$ лежит по одной
целой точке, откуда $m+m^*=4+8=12$.

\bigskip

$$
\epsfysize=3cm
\epsfbox{7b.eps}
$$

\nopagebreak
\smallskip

\centerline{Рисунок~2.}

\bigskip

Пусть теперь $M=A_1\dots A_n$.
Будем считать все целые точки на сторонах многоугольника $M$
вершинами (возможно, с углом $180^\circ$).
На определении $M^*$ это не отразится.
Предположим, что некоторый треугольник $A_{i-1}A_iA_{i+1}$
{\it прост}, то есть не содержит целых точек, отличных от вершин
(ни внутри, ни на сторонах). {\it Элементарной операцией}
$A_1\dots A_{i-1}A_iA_{i+1}\dots A_n \to A_1\dots A_{i-1}A_{i+1}\dots A_n$
назовем операцию удаления треугольника $A_{i-1}A_iA_{i+1}$
из $M$ или обратную к ней.
Наш метод сведения основан на следующем утверждении:

\smallskip
(2) {\it При элементарной операции величина $m+m^*$ сохраняется.}
\smallskip

Нам достаточно доказать, что, например, удаление простого
треугольника $A_1A_2A_3$ из $M$ приводит к добавлению простого
треугольника $A_{12}A_{13}A_{23}$ к $M^*$ (Рис.~3).
Здесь через $A_{kl}$ обозначена точка, такая что
$\overrightarrow{OA_{kl}}=\overrightarrow{A_kA_l}$.
В частности, если $l=k+1$, то $A_{kl}$ --- вершина многоугольника $M^*$.
Удалим $A_1A_2A_3$. Тогда у многоугольника $M^*$ исчезнут
вершины $A_{12}$ и $A_{23}$, зато добавится новая вершина $A_{13}$.
Ее еще нужно соединить отрезками с $A_{n1}$ и $A_{34}$.
Покажем, что точки $A_{12}$ и $A_{23}$ лежат на этих отрезках.
В самом деле, так как $O$ --- единственная целая точка внутри $M$,
то треугольники $A_1OA_3$, $A_2OA_3$, $A_4OA_3$ --- простые.
Из формулы Пика следует, что их площади равны $1/2$.
Поскольку они имеют общее основание $OA_3$, то проекции
векторов $\overrightarrow{A_1A_3}$,
$\overrightarrow{A_1A_3}$ и $\overrightarrow{A_1A_3}$
на перпендикуляр к $OA_3$ равны.
Отсюда следует, что точки $A_{13}$, $A_{23}$ и $A_{34}$ лежат
на одной прямой, а так как $M$ --- выпуклый, то
$A_{23}$ лежит между двумя остальными. Аналогично доказывается,
что $A_{12}$ принадлежит отрезку $A_{n1}A_{13}$.
Значит, преобразование $M^*$ сводится к добавлению треугольника
$A_{12}A_{13}A_{23}$.
Заметим, что треугольник $OA_{12}A_{13}$ получается из
простого треугольника $A_1A_2A_3$ параллельным переносом,
а $OA_{23}A_{13}$ --- центральной симметрией. Поэтому
треугольник $A_{12}A_{13}A_{23}$ --- простой, что и требовалось.

\bigskip

$$
\epsfysize=3cm
\epsfbox{9.eps}
$$

\nopagebreak
\smallskip

\centerline{Рисунок~3.}

\bigskip


Для доказательства теоремы остается заметить следующее:

\smallskip
(3) {\it Из любого многоугольника $M$ серией элементарных операций
можно получить параллелограмм без целых точек на сторонах.}
\smallskip

Действительно, предположим вначале, что у $M$ есть диагональ,
не проходящая через $O$. Разрежем $M$ вдоль этой диагонали и
рассмотрим ту из полученных частей, которая не содержит $O$.
Эта часть обязательно содержит простой треугольник вида $A_{i-1}A_iA_{i+1}$.
Поэтому элементарной операцией можно уменьшить число $m$.
Будем действовать так, пока это возможно.
Очевидно, есть только три случая, когда требуемой диагонали не найдется:

\noindent A) $m=4$, $M=ABCD$, $O=AC\cap BD$. Так как отрезки
$OA$, $OB$, $OC$ и $OD$ не содержат целых точек, то $OA=OC$ и $OB=OD$,
то есть $ABCD$ --- искомый параллелограмм.

\noindent B) $m=4$, $M=ABCD$, один из углов, скажем, $BCD$ --- развернутый.
В этом случае обозначим через $D'$ точку, симметричную
точке $D$ относительно $O$, через $E$ --- середину $D'B$.
Искомая серия элементарных операций имеет вид:
$$
ABCD\to AEBCD\to AD'EBCD\to AD'ECD\to AD'CD \quad\text{(Рис.~ 4).}
$$

\noindent C) $m=3$, $M=ABC$. В этом случае обозначим через $A'$ и $C'$
точки, симметричные относительно $O$ вершинам $A$ и $C$ соответственно.
Тогда искомая серия имеет вид:
$$
ABC\to AC'BC\to AC'BA'C\to AC'A'C \quad\text{(Рис.~ 4).}
$$
Теорема доказана.

\bigskip

$$
\epsfysize=3cm
\epsfbox{10.eps}
$$

\nopagebreak
\smallskip

\centerline{Рисунок~4.}

\bigskip

{\bf Благодарности.} Авторы благодарны В.В. Прасолову
за постановку задачи.